\documentclass[12pt,letterpaper, reqno]{amsart}

\usepackage{amssymb,latexsym, amsmath, amsxtra}
\usepackage[dvips]{graphics}

\theoremstyle{plain}
\newtheorem{theorem}{Theorem}[section]
\newtheorem{corollary}[theorem]{Corollary}

\newtheorem{principle}[theorem]{Principle}

\newtheorem{lemma}[theorem]{Lemma}

\theoremstyle{definition}

\theoremstyle{remark}

\numberwithin{equation}{section}
\numberwithin{theorem}{section}
\numberwithin{table}{section}
\numberwithin{figure}{section}

\newcommand{\C}{\mathbb C}

\DeclareMathOperator{\sinc}{sinc}

\newcommand{\term}[1]{\textbf{#1}}

\def\({\left(}
\def\){\right)}

\begin{document}
\title[Jensen polynomials and RH]{Jensen polynomials are not a plausible route \\ to proving the Riemann Hypothesis}
\author{David~W.~Farmer}
\address{American Institute of Mathematics,
600 East Brokaw Road,
San Jose, CA 95112-1006}
\email{farmer@aimath.org}


\begin{abstract}
Recent work on the Jensen polynomials of the Riemann xi-function
and its derivatives found a connection to the Hermite polynomials.
Those results have been suggested to give evidence for the
Riemann Hypothesis, and furthermore it has been suggested that those
results shed light on the random matrix statistics for zeros of
the  zeta-function.  We place that work in the context of
prior results, and explain why the appearance of Hermite polynomials
is interesting and surprising, and may represent a new type of
universal law which refines M.~Berry's ``cosine as a universal attractor''
principle.  However, we find there is no justification
for the suggested connection to the Riemann Hypothesis, nor for the
suggested connection to the conjectured random matrix statistics for
zeros of L-functions.  These considerations suggest that Jensen polynomials,
as well as a large class of related polynomials, are not useful for
attacking the Riemann Hypothesis.
We propose general criteria for determining whether an equivalence to the
Riemann Hypothesis is likely to be useful.
\end{abstract}

\thanks{This research was supported by the National Science Foundation.}


\maketitle

\section{Introduction} 

Two recent papers \cite{GORZ, Bo} revisit the classical result of
Jensen~\cite{J, P} that the Riemann Hypothesis(RH) is true if and only if
all of the associated Jensen polynomials, defined in \eqref{eqn:jensencl}
below, have only real zeros.  The two recent papers actually concern
another version of the Jensen polynomials, which we call the ``even''
Jensen polynomials, defined in \eqref{eqn:jenseneven}.  An interesting
connection was found with the Hermite polynomials.

In this paper we examine the recent work on Jensen polynomials
in the context of prior work on repeated
differentiation of entire functions~\cite{Ber, FR, Ki, Kim},
and on differentiation-like operations~\cite{Do, RT}.
That perspective
explains why the new connection to Hermite polynomials is
interesting, but it also suggests why there is no connection to the
Riemann Hypothesis nor to the random matrix statistics of zeros of
the zeta function.
These considerations further suggest that the Jensen polynomials,
as well as a large class of related polynomials, are not a useful tool for
approaching the Riemann Hypothesis. 
We introduce terminology
which can serve as a guide to deciding whether an equivalence to RH
is likely to be useful for resolving the Riemann Hypothesis, or if
the equivalence is just a curiosity.

\section{The classical Jensen polynomials}\label{sec:jensencl}
Suppose
\begin{equation}\label{classicaltaylor}
f(z) = \sum_{j=0}^\infty \frac{\alpha(j)}{j!} z^j
\end{equation}
is an entire function of order less than two. 
One can associate the
\term{$d$th classical Jensen polynomial for the $n$th derivative of $f$},
given by
\begin{equation}\label{eqn:jensencl}
J^{d,n}_{f,cl}(z) := \sum_{j=0}^d \binom{d}{j} \alpha(j + n) z^j .
\end{equation}
The ``$cl$'' in the subscript refers to these polynomials being
``classical'' in the sense that \eqref{eqn:jensencl} is the standard
definition of the Jensen polynomials.  An alternate notation for
those polynomials is $J^{d,n}_{\alpha,cl}$, where the first subscript
refers to the Taylor series coefficients instead of to the function.
We will also consider the ``even'' Jensen polynomials, defined
in~\eqref{eqn:gamma}.

One reason for interest in the classical Jensen polynomials is:
\begin{enumerate}
\item \label{jensenlim} $\displaystyle \lim_{d \to \infty} J^{d,n}_{f,cl}(z/d) = f^{(n)}(z)$,
with uniform convergence for $z$ in a compact set, and
\item \label{rhandjensen} $f$ has only real zeros if and only if $J^{d,0}_{f,cl}$ 
has only real zeros for
all $d$.
\end{enumerate}
Note that item (\ref{jensenlim}) directly gives one of the implications
in item~(\ref{rhandjensen}).
For real entire functions of order less than two,
the property of having only real zeros is preserved under differentiation,
so an equivalent reformulation of item~(\ref{rhandjensen}) is that
$f$ has only real zeros if and only if $J^{d,n}_{f,cl}$ 
has only real zeros for all $d$ and all~$n$.

We will describe results in the literature as they apply to the classical
Jensen polynomials $J^{d,n}_{f,cl}$ as $n \to \infty$, and then consider
the corresponding problem for the even Jensen polynomials
considered in~\cite{GORZ, Bo}.

For the functions under consideration here,
differentiation preserves real zeros.  Much more is true.  A beautiful
result of Kim~\cite{Kim} asserts that if $f$ is an entire
function of order less than~$2$, which is real on the real axis, and which
has all zeros in a strip $|\Im(z)| < A$, then for any fixed $R>0$, 
if $n$ is sufficiently large then $f^{(n)}$ has only real zeros in $|z| < R$.
In other words, in any compact region, if you differentiate such functions
enough times, all zeros are real.  A corollary is
that for any $d$, if $n$ is large enough then the classical Jensen polynomial
$J^{d,n}_{f,cl}$ has only real zeros.

For a large subset of the functions for which Kim's theorem applies,
even more is conjectured:  not only do the zeros move to the real axis,
they also approach equal spacing.  Since (up to a simple change of variables)
the only even, real, entire function of order less than~$2$ with equal spaced real
zeros is the
cosine function, it is conjectured that for a large class of functions,
repeated differentiation leads to the cosine function, up to a simple rescaling.
A precise form of this conjecture was made by Berry~\cite{Ber}, who phrased it
as 
\begin{center}
\emph{
$\cos(\omega_n t + \delta_n)$ is a universal attractor of the derivative map
},
\end{center}
and by
Farmer and Rhoades~\cite{FR} from a slightly different perspective
based on the density of zeros of the function.

A relevant instance of that conjecture was proven by Ki~\cite{Ki}.
Let
$$
\Xi(z) = \xi(\tfrac12 + i z)
$$
be the Riemann $\Xi$-function, where
$$
\xi(s) = {\textstyle \frac12} s (1-s)\pi^{-s/2} \Gamma\left(\frac{s}{2}\right) \zeta(s).
$$
The function $\Xi$ is even and is real on the real axis,
and has all zeros in the strip $-\frac12 < \Im(z) < \frac12$.  Thus, it is a theorem
that all zeros of $\Xi^{(n)}(z)$ for  $|z| < T$ are real if $n$ is sufficiently large,
and so for each $d$, if $n$ is large enough the classical Jensen polynomial
$J^{d,n}_{\Xi,cl}$ has only real zeros.
It was further conjectured~\cite{FR} that, suitably rescaled, $\Xi^{(n)}(z)$
approaches $\cos(z)$.  That conjecture was proven by Ki~\cite{Ki}:
\begin{theorem}[Ki \cite{Ki}]\label{thm:ki}
There exist positive decreasing sequences $A_n$ and $C_n$ such that
\begin{equation}
\lim_{n\to\infty} (-1)^n A_n \,\Xi^{(2n)}(C_n z) = \cos(z),
\end{equation}
uniformly on compact subsets of~$\C$.
\end{theorem}
That theorem also follows from a proposition of Coffey~\cite{C}.
The analogous result holds for functions in the extended Selberg class~\cite{GH}.
Functions in that class have a functional equation but  not necessarily
an Euler product, and so it includes many examples that do not satisfy
the analogue of the Riemann Hypothesis.  The same result holds for
random functions~\cite{PS}, which by construction satisfy the analogue of
the Riemann Hypothesis but have Poisson statistics for their zeros.
Conrey's result~\cite{Co} that $\Xi^{(n)}$ has $(100 - O(1/n^2))$ percent of
its zeros on the real axis is a quantitative version of Kim's theorem,
and can be seen as a foreshadowing of Ki's result.

Theorem~\ref{thm:ki} implies that the rescaled Taylor series coefficients
of~$\Xi^{(n)}$ converge to those of cosine.
Since the Taylor coefficients of cosine have a simple form,
the Jensen polynomials of cosine can be written explicitly:
\begin{equation}\label{eqn:jensencosine}
J^{d,0}_{\cos,cl}(z) = \frac12 ((1+ i z)^d + (1 - i z)^d).
\end{equation}
By Theorem~\ref{thm:ki} and \eqref{eqn:jensencosine} we have
\begin{corollary}
There exist positive decreasing sequences $A_n$ and $C_n$ such that
$$
\lim_{n\to\infty} (-1)^n A_n J^{d,2n}_{\Xi,cl}(C_n z) = \frac{(1+ i z)^d + (1 - i z)^d}{2}
$$
as $n\to \infty$.  In particular, for each $d$, if $n$ is sufficiently large then $J^{d,2n}_{\Xi,cl}$
has only real zeros.
\end{corollary}

We see that the classical Jensen polynomials $J^{d,n}_{f,cl}$
having real zeros
for large $n$ is a general phenomenon, following from the fact that,
for a large class of entire functions, repeated differentiation leads
to the cosine function.
In particular, differentiation causes a loss of information about
the zeros of the functions considered here, and so in terms of the
Riemann Hypothesis there is little revealed
by the derivatives of the function.
In Section~\ref{sec:partialRH}
we elaborate with an illustrative example.  But first we consider
a different form of the Jensen polynomials.

\section{The even Jensen polynomials}
If $f$ is an even function, it is natural to write
\begin{equation}\label{eqn:gamma}
f(z) = \sum_{j=0}^\infty \gamma(j) \frac{z^{2j}}{j!}  .
\end{equation}
From this we define the \textbf{even Jensen polynomials}, which are the
subject of~\cite{GORZ, Bo}:
\begin{equation}\label{eqn:jenseneven}
  J^{d,n}_{f,ev}(z) := \sum_{j=0}^d \binom{d}{j} \gamma(j+n) z^j.
\end{equation}
As in the classical case, the first subscript could be the even Taylor coefficients,
$\gamma$, instead of the function.

Note that the even Jensen polynomial of $f(z)$ is the classical Jensen
polynomial of $f(\sqrt{z})$.
In the case of the Riemann $\xi$-function, the Riemann hypothesis is
equivalent to the assertion that $\xi(\frac12 + \sqrt{z})$ has zeros
only on the negative real axis, or equivalently, $\Xi(\sqrt{z})$
has zeros only on the positive real axis.

The terminology of ``classical'' and ``even'' Jensen polynomials is not
standard, but we felt the terminology was necessary in order to avoid confusion. 
We write $J^{d,n}$ when we wish to make a statement that applies
to either case.

The main results of~\cite{GORZ} are precise asymptotics for~$\xi^{(2n)}(\frac12)$
and a new phenomenon relating asymptotic properties of certain sequences
to the Hermite polynomials.  Those results combine to produce:
\begin{theorem}[Griffin, Ono, Rolen, and Zagier \cite{GORZ}] There exist sequences $\mathcal A_n$, $\mathcal B_n$, and $\mathcal C_n$
such that
$$
\lim_{n\to\infty} \mathcal A_n J^{d,n}_{\xi,ev}(\mathcal C_n z + \mathcal B_n) = H_d(z) ,
$$
uniformly for $z$ in a compact subset of $\C$,
where $H_d$ is the $d$th Hermite polynomial and the
subscript $\xi$ refers to $\xi(\frac12 + z)$.
\end{theorem}

We compare this to Ki's theorem quoted above, which implies the following:
\begin{corollary}  
There exist sequences $A_n$ and $C_n$ such that
$$
\lim_{n\to\infty} A_n J^{d,n}_{\xi,ev}(C_n z) = (1 + z)^d ,
$$
uniformly for $z$ in a compact subset of $\C$.
\end{corollary}

\begin{figure}[htp]
\scalebox{0.6}[0.6]{\includegraphics{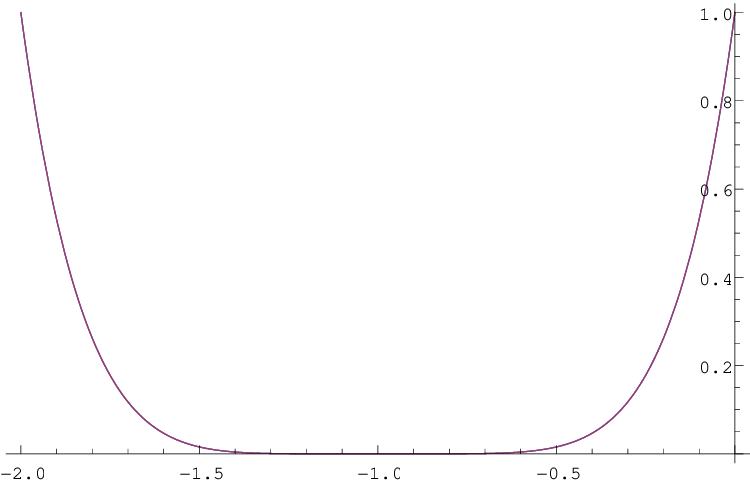}}
\hskip 0.4in
\scalebox{0.6}[0.6]{\includegraphics{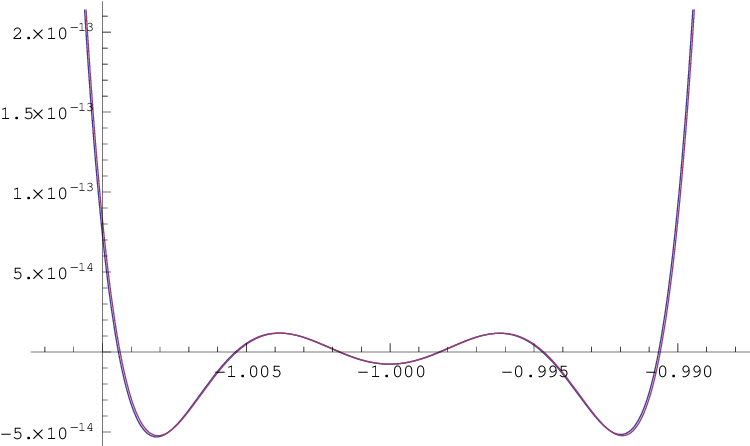}}
\caption{\sf
The even Jensen polynomial $J^{6,10000}_{\xi,ev}(x)$,
rescaled as described in the text,
for $-2\le x \le 0$ on the left, and $-1.012 \le x \le -0.988$
on the right.
The plot on the left is superimposed with the graph of $(1 + x)^6$,
and the plot on the right is superimposed with the Hermite polynomial
$H_6(x)$, shifted and
scaled.
} \label{fig:J6of10k}
\end{figure}

How can we reconcile the fact that the $d$th even Jensen polynomials simultaneously
converge both to $(1 + z)^d$ and to $H_d(z)$?  This apparent
conundrum is easily resolved by examining a plot of
the polynomial.  In Figure~\ref{fig:J6of10k} we show graphs
of $A_n  J^{d,n}_{\xi,ev}(C_n x)$ for $d=6$, $n=10000$. That is,
the 6th even Jensen polynomial of the 10000th derivative of $\xi(\frac12 + z)$.
The plot on the left covers the range $-2 \le x \le 0$, and the plot
on the right covers $-1.012 \le x \le -0.988$.  

Each plot in Figure~\ref{fig:J6of10k}  actually
contains a superposition with a second graph:  $(1+x)^6$ on the left,
and $H_6(x)$, shifted and scaled, on the right.  In both cases the plots are so close that
the two graphs are indistinguishable to the eye.
 
We see that the main result of~\cite{GORZ} contains more information than
the theorem of Ki~\cite{Ki} because $H_d(x)$, suitable shifted and scaled,
looks just like $(1+x)^d$, but the converse is not true.

We suggest that the results in~\cite{GORZ} can be interpreted as a refinement
of the general ``cosine universality'' of Berry and Farmer-Rhoades.  That is:
\begin{principle}[``Hermite Universality'']\label{prin:HU}
For a large class of functions, not only does repeated differentiation
lead to the (rescaled) cosine function, but the convergence occurs
in a particularly regular and uniform way, characterized by the
appearance of the Hermite polynomials within the shifted
and rescaled even Jensen polynomials.
\end{principle}

We can be more specific about what this principle predicts.  Suppose $f(z)$
is an even real entire function for which Cosine Universality should
hold.  Interpreting Cosine Universality as a statement about
Taylor coefficients,  we see that (suitably scaled but not shifted),
the $n$th derivative of $f(\sqrt{z})$ approaches $e^{-z}$, and so 
$J^{d,n}_{f, ev}(z)$
(suitably scaled but not shifted) approaches $(1-z)^d$ as $n\to\infty$.
Taken at face value, that limit does not directly imply that
$J^{d,n}_{f, ev}(z)$ has only real zeros for sufficiently large~$n$ 
(although one might conclude that from other considerations).
Hermite Universality does imply that $J^{d,n}_{f, ev}(z)$ has
only real zeros for sufficiently large~$n$, and it further implies that
those zeros are arranged like the zeros of a Hermite polynomial,
shifted and scaled into a small interval around $z=1$.

Griffin, Ono, Rolen, and Zagier~\cite{GORZ} verify this principle in many cases,
and also consider it as applied to sequences that are not
being viewed  as the derivatives of an entire function.

Note that the principle is not restricted to even functions.
However the concept of ``even'' Jensen polynomial has yet to be
defined for functions which are not even, but which when repeatedly
differentiated and slightly shifted, converge to the cosine function.
Presumably there are functions for which Cosine Universality applies
but Hermite Universality does not -- perhaps functions without sufficient
regularity in the spacings of their zeros.

There are a couple of facts that point to Principle~\ref{prin:HU} as an interpretation of the
results in~\cite{GORZ}.  First is that the work of Ki~\cite{Ki},
its generalization to the extended Selberg class~\cite{GH}, 
the proof of Newman's conjecture and its generalization~\cite{Do, RT},
and the
work under discussion~\cite{GORZ}, all rely on the fact that
functions under consideration can be written in a form
similar to
\begin{equation}\label{eqn:Xiexp} 
\Xi(z) = \int_{-\infty}^\infty \varphi(u) e^{i z u} \, du
\end{equation}
where $\varphi$ decreases rapidly.  Such an expression is amenable
to analyzing derivatives of $\Xi$, and~\cite{GORZ} carries the
analysis farther than previous efforts.

The second reason comes from considering some simple examples,
which we describe after initial preparations in the next section.

\section{Not all equivalences to RH are created equal}\label{sec:partialRH}

We have seen that the Jensen polynomials of derivatives, $J^{d,n}$ for
$n\ge 1$, do not shed any light on the Riemann Hypothesis, because each
increase in the differentiation index loses information about the
location of the zeros.  In this section we give
another example to further illustrate that point,
but our main purpose now is
to complete the claim in the title of this paper, describing why $J^{d,0}$,
with differentiation index 0, is also not
a useful tool for exploring RH.

Our argument is in three parts.  First we divide the equivalences to RH into
different categories.  Then we suggest criteria for deciding, within each category,
whether an equivalence is likely to be helpful for proving~RH.  In particular, we make the point that
some equivalences to~RH are unlikely to be useful for proving~RH.

Given this perspective, we then consider the case of Jensen polynomials and
a family of related equivalences.

\subsection{Towards a taxonomy of equivalences}\label{sec:taxonomy}
Many equivalences to RH fall into one or more of the following categories.

\begin{enumerate}
\item[(A)] A subset or superset of an existing equivalence.
(In the case of superset, there are two subcategories, depending on whether
or not
the additional conditions are logical consequences of the previous conditions.)

\item[(B)] A repackaging of an existing equivalence.

\item[(C)] A translation into a different language.
\end{enumerate}

In case (A), it is reasonable to interpret a subset equivalence as a promising
route to proving RH, because there are fewer conditions to satisfy.  And a superset
equivalence, if the additional conditions are not logical consequences of the
existing conditions, can be interpreted as a promising route to disproving~RH,
because there are more opportunities to obtain a contradiction.
Thus, except in the case where simple logic indicates that the extra conditions
provide no additional information, such equivalences cannot be easily ruled out
as a plausible route to resolving~RH.

In case (B), the potential usefulness of the equivalence hinges on
whether the information in the previous equivalence has been concentrated
or dispersed.  We illustrate the idea with the 
equivalences of Robin~\cite{Rob} and Lagarias~\cite{Lag}.
Those equivalences involve upper bounds 
of the form
\begin{equation}\label{eqn:robinlagarias}
\sigma(n) \le f(n)
\end{equation}
where
$\sigma(n) = \sum_{d|n} d$ is the divisor sum function and
$f$ is given explicitly.  The proofs of those
equivalences start with the RH equivalence involving the error term in
the prime number theorem:
\begin{equation}\label{eqn:pili}
\pi(x) = \mathrm{Li}(x) + O(x^{\frac12 + \varepsilon}).
\end{equation}
A violation of~\eqref{eqn:pili} for a particular $x$ is used to exhibit an
integer $n$ where $\sigma(n)$ is particularly large:
large enough to violate~\eqref{eqn:robinlagarias}.  The relevance to
our discussion here is that \emph{the integer $n$ is enormously larger than~$x$}.
Thus we say the equivalence has \term{dispersed} the information:
the new condition requires searching further in order to obtain the same
information which was previously available.
The dispersal of information is an indication that the equivalence is unlikely to be
helpful for resolving~RH.

In case (C), the issue is whether the translation could allow
the use of new tools.
An example which \emph{does} afford new tools is the equivalence between RH and~\eqref{eqn:pili}.
Indeed, the Prime Number Theorem is equivalent to the nonvanishing of
the $\zeta$-function on the line $\sigma=1$, and both parts of the equivalence
have been proven independently.

One could view Robin's and Lagarias' equivalences
as falling into case (C), since $\sigma(n)$ does not literally appear in~\eqref{eqn:pili}.
However, the use of $\sigma(n)$ is just convenient packaging, and there are no
special properties of the $\sigma$-function which are relevant to the proof.

We consider one more example of case (C) before
returning to the  Jensen polynomials.
\begin{lemma}\label{farmerequivalence} The following are equivalent:
\begin{enumerate}
\item\label{rhsz} The Riemann Hypothesis is true and all zeros of the $\zeta$-function are simple,

\item\label{taylorz}
For all $R > 0$, if $n \ge n(R)$ is sufficiently large, then
inside the disc $|z| < R$ the $n$th order Taylor polynomial
for $\Xi(z)$ has only real zeros.
\end{enumerate}
\end{lemma}

\begin{proof}
If (\ref{rhsz}) is true, then (\ref{taylorz}) follows from Taylor's theorem,
Rouch\'e's theorem, and the fact that the Taylor coefficients of the $\Xi$-function are real.

In the other direction, suppose the $\Xi$-function had a multiple zero.  By
Theorem~\ref{thm:ki}, for large $n$ the signs of $\Xi^{(2n)}(0)$ alternate,
so in a neighborhood of the multiple zero the $(2n)$th order Taylor polynomial
is alternately larger and smaller than
the $\Xi$-function.  So a double zero of $\Xi$ would alternately be a pair of real zeros
and a pair of complex zeros of its $(2n)$th order Taylor polynomial,
and a higher odd-order zero would only contribute a single
real zero to the Taylor polynomial.
\end{proof}

Does that equivalence to (RH + simple zeros) open the possibility of applying new tools to the
problem?  The answer might not be as definitive as in the previous examples,
but (in the author's opinion) it seems fairly clear that nothing has been
gained by translating to Taylor polynomials.

Thus, in each of cases (A), (B), and (C) we have criteria to judge whether or not
a given equivalence is a plausible route to resolving~RH.  We do not claim 
to ``prove'' that an equivalence cannot be used to resolve~RH,
but mathematics is a human endeavor, and human effort is limited, so it is
helpful to have reasons for deciding what effort is likely to be fruitful.
A similar sentiment was expressed by Poincar\'e more than 100 years ago~\cite{Po}:

\begin{quote}
For a construction to be useful and not mere waste of mental effort,
for it to serve as a stepping-stone to higher things, it must first of all
possess a kind of unity enabling us to see something more than the
juxtaposition of its elements.
\end{quote}

In the next section we explain why the Jensen polynomials are even less
useful than the Taylor polynomials as an approach to~RH.

\subsection{Jensen polynomials disperse the information}
There is evidence in the literature that $J^{d,0}$
is not effective at detecting violations of the Riemann Hypothesis.  Namely,
Chasse~\cite{Ch} proved that if all the zeros $\rho = \beta + i \gamma$
of the zeta-function are on the critical line for $|\gamma |<T$,
then $J^{d,0}$ has only real zeros for $d < T^2$.  
In other words, Jensen polynomials disperse the
information about zeros.  We illustrate this idea with a simple example.

\begin{figure}[htp]
\includegraphics{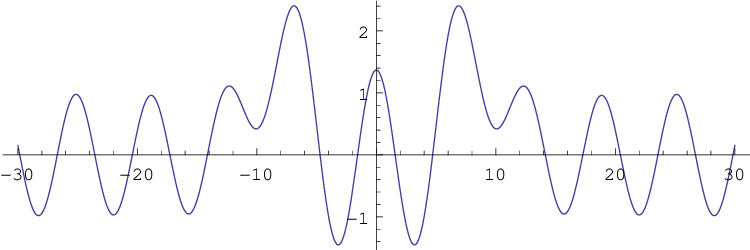}
\caption{\sf
A graph of the function obtained by moving the zeros
at
$x = \pm \frac{5\pi}{2}$ and $\pm \frac{7\pi}{2}$ of $\cos(x)$
to $\pm 10 \pm i$.
} \label{fig:X10}
\end{figure}

Consider a function which is entire of order 1, even, real on the real axis,
and has all its zeros in a strip $-A < \Im(z) < A$.  The analogue of the
Riemann Hypothesis is that all of the zeros are real.  An example,
which presumably has only real zeros, is
the Riemann $\Xi$-function.  
An example which does not 
have only real zeros is
$$
X_{10}(z) = \cos(z)
   \frac{(z^2-(10+i)^2)(z^2-(10-i)^2)}{(z^2-(\frac{5\pi}{2})^2)(z^2-(\frac{7\pi}{2})^2)}.
$$
In words, $X_{j}(z)$ is the function obtained when the pairs of zeros of
$\cos(z)$ closest to $\pm j$ are moved to $\pm j \pm i$, and above is a formula for
$X_{10}$.  Figure~\ref{fig:X10}
shows a graph of $X_{10}(x)$.

Examining the graph of $X_{10}$, it can be seen that all zeros
of the first derivative, $X_{10}'$, are real, therefore the same is true of all higher
derivatives.  Thus, the classical Jensen polynomial
$J^{d,n}_{X_{10}, cl}$ has only real zeros for all $n\ge 1$,
as does the even Jensen polynomial $J^{d,n}_{X_{10}, ev}$
for all even $n\ge 2$.

But what about $J^{d,0}_{X_{10}, cl}$? We know that this will have
non-real zeros if $d$ is large enough, but how large is large enough?
The first two zeros (in magnitude) of $X_{10}$ are real, and then a pair of
complex conjugate zeros.  Since the ``Riemann Hypothesis''
fails for $X_{10}$ almost immediately,
one might guess
that $J^{d,0}_{X_{10}, cl}$
should have a non-real zero for $d$ quite small. This is not
the case.  By a direct calculation (we used Mathematica),
$J^{d,0}_{X_{10}, cl}$ has only real zeros for~$d\le 118$,
and for all larger $d$ it has non-real zeros.  

Table~\ref{tab:Xjhyperbolic} shows, for various $X_j$, the maximal $d$
such that $J^{d,0}_{X_j, cl}$ has only real zeros.
The data in that table confirm the impression from Chasse's theorem,
that Jensen polynomials are inefficient at detecting non-real zeros.
Furthermore, the Jensen polynomials, which are defined in terms of the
Taylor series coefficients, are not efficient at extracting information
from those coefficients.  The $d$th order Taylor polynomial of
$X_j$ also detects the non-real zero if $d$ is large enough,
in accord with Lemma~\ref{farmerequivalence}:  this is shown
in the bottom row of Table~\ref{tab:Xjhyperbolic}.  We see that
the Jensen polynomials require significantly more Taylor coefficients than
the Taylor polynomials to detect the non-real zeros.

\begin{table}[]
\begin{tabular}{lllll}
\ \ \ \ \ \ \ \ \ \ \ \ \ \ \ \ \ \ $j$  & 10  & 20 & 40 & 60 \\
\hline
\# first zeros of $X_j$ are real         & 2   & 4  &  12  &  18  \\
$J^{d,0}_{X_j, cl}$ has only real zeros for $d\le$    & 118 & 749   & 1897   &  4242  \\
$d$th  Taylor  polynomial detects non-real zero for $d \ge$ \ \ \ \  & 20  &  60  & 118  &  175 \\
 & & & & \\
\end{tabular}
\caption{Tabulating the relative effectiveness of the Jensen
polynomials and the Taylor polynomials for detecting violations
of the Riemann Hypothesis, using the function $X_j$ as a model.}\label{tab:Xjhyperbolic}
\end{table}

In the terminology of Section~\ref{sec:taxonomy}, Jensen polynomials are a repackaging of
the equivalence in Lemma~\ref{farmerequivalence}.  And since the Jensen polynomials disperse
the information in the Taylor polynomials, we are justified in asserting that
the Jensen polynomials are even less useful than the Taylor polynomials
as a tool for resolving~RH.

A similarity between Jensen and Taylor polynomials is that they approximate
$\Xi(z)$ when $|z|$ is small.  
It is tempting to view the Jensen polynomials as ``better'' because for larger $z$
the zeros of Jensen polynomials are real, while Taylor polynomials tend to have
many complex zeros.  But, that apparently nice property is
just a distraction.  
One set of functions has meaningless zeros on a line, and the other has
meaningless zeros near a circle.
In both cases the extraneous zeros say very little about the
function being approximated.  The apparently nice property of having extra real
zeros comes at the cost of converging to the function more slowly.  One must distinguish
between elegance in the statement of a proposition, and actually being useful as a tool
to prove new results.  That criticism also applies to the equivalences due to
Robin and to Lagarias.

\subsection{Other Jensen-like polynomials}

There are other polynomials generated from the Taylor coefficients
which only have real zeros if and only if the original function has
only real zeros.  For example, a recent paper of O'Sullivan~\cite{OS}
considers the polynomials
\begin{equation}\label{eqn:cormacpoly}
   P^{d,n}(z) := \sum_{j=0}^d \binom{d}{j} \gamma(j + n) H_{d-j}(z).
\end{equation}
In other words, the Jensen polynomial with $z^j$ replaced by the
Hermite polynomial~$H_j(z)$.  

O'Sullivan shows that 
these polynomials have the same property that make the Jensen polynomials
interesting: 
$\Xi(z)$ has only real zeros if and only if
$P^{d,n}$ has only real zeros for all $d$,~$n$.  That result is a special
case of a more general result whereby any element of the Laguerre-P\'olya
class produces a sequence of polynomials which can be put in place
of the Hermite polynomials in~\eqref{eqn:cormacpoly}.  Thus, there
is a wealth of seemingly different sequences of polynomials,
any one of which can detect a violation of the Riemann Hypothesis.

We have argued that the Jensen polynomials are not a useful tool
for attacking the Riemann Hypothesis.  Might one of those other
sequences of polynomials turn out to be more useful?
Sadly,~no.  O'Sullivan goes on to show that if the even Jensen polynomial
$J^{d,n}_{\Xi,ev}$ has only real zeros, then so does~$P^{d,n}$.
In other words, $P^{d,n}$ is less useful at detecting violations
of the Riemann Hypothesis.  The proof in~\cite{OS}
is in the context of~$P^{d,n}$, but presumably the analysis
extends to all the other sequences of polynomials.

\subsection{Other differentiation-like operations}\label{sec:newman}
de Bruijn~\cite{dB} and Newman~\cite{N} considered the following operation, which
uses the notation of \eqref{eqn:Xiexp},
\begin{equation}\label{eqn:Xiexpt} 
\Xi_t(z) = \int_{-\infty}^\infty e^{t u^2} \varphi(u) e^{i z u} \, du .
\end{equation}
The de Bruijn-Newman constant is defined by
$\Lambda = \inf\{t\ :\ \Xi_t \text{ has only real zeros}\}$.
Since $\Xi_0 = \Xi$, the Riemann Hypothesis is equivalent to $\Lambda \le 0$.
Newman conjectured $\Lambda \ge 0$, which was proven recently by
Rodgers and Tao~\cite{RT}.  That result was generalized to
the extended Selberg class (most of which does not satisfy the
analogue of the Riemann Hypothesis) by Dobner~\cite{Do}.  
As Dobner notes, the method ``does not require any information about the zeros''
of the function.  In particular, these results say nothing about 
Lehmer pairs of zeros, which is somewhat ironic since previously the
lower bounds on $\Lambda$ came from Lehmer pairs.

The de~Bruijn-Newman operation
$\Xi \to \Xi_t$ has several properties in common with differentiation
$\Xi \to \Xi^{(j)}$.  For example, if $\Xi_{t_0}$ has only real zeros
then $\Xi_t$ has only real zeros for all~$t > t_0$.  Also, if $t > 0$
then as $x\to \infty$ the zeros of $\Xi_t(x)$ approach equal spacing.
Thus, the de~Bruijn-Newman operation is even more efficient
than differentiation at losing information and causing the zeros to
approach equal spacing.

In some sense, de~Bruijn-Newman operation is like
repeated differentiation $\Xi^{(j)}(z)$ where $j$ is an increasing
function of~$z$.  It is possible to be somewhat precise about that
remark.  As shown in Section~3 of~\cite{Do}, if $z$ is real then
the main contribution in~\eqref{eqn:Xiexpt}
is concentrated near $u \approx z$. If $n$ is the integer closest to $t z^2$,
then the largest term in the Taylor series for $e^{t u^2}$ is approximately
\begin{equation}
\frac{t^n}{n!} u^{2n} \approx \frac{t^n}{n!} u^{2 t z^2} .
\end{equation}
In other words, as $z\to\infty$ the zeros of $\Xi_t(z)$ are approaching
equal spacing are a rate comparable to the $2 t z^2$rd derivative
$\Xi^{(2 t z^2)}(z)$.  That analysis
may not be rigorous, but it does explain why
the de~Bruijn-Newman operation is extremely effective at causing
zeros to become equally spaced.

The above discussion was intended to emphasize the point that functions
with a representation similar to \eqref{eqn:Xiexp} have a nice distribution
to their zeros, which becomes nicer under operations similar to differentiation.
But those seemingly magical properties are from the realm of analysis,
not number theory, and those properties hold whether or not
the function satisfies a Riemann Hypothesis.

\section{Hermite polynomials and the function $X_{10}$}

The function $X_{10}$ in the previous section approaches cosine under
repeated differentiation.  We now view that function through the
lens of the results in~\cite{GORZ}.
Let $X_{10}(x) = \sum \alpha(n) x^n / n!$.  Using Cauchy's theorem,
as described in~\cite{Bor},
we computed the $\alpha(n)$ for a few $n$ near $100000$, shown in
Table~\ref{tab:alphan}.
\begin{table}[]
\begin{tabular}{cc}
$n$         & $\alpha(n)$ \\
\hline
100000 & $\phantom{-} 1.0000000015411856213980266026829$ \\
100002 & $-1.0000000015411239767472129801496$ \\
100004 & $\phantom{-} 1.0000000015410623357948380052808$ \\
100006 & $-1.0000000015410006985406058281322$ \\
100008 & $\phantom{-} 1.0000000015409390649842206283415$ \\
100010 & $-1.0000000015408774351253866151245$ \\
100012 & $\phantom{-} 1.0000000015408158089638080272719$ \\
 & \\
\end{tabular}
\caption{The nonvanishing Taylor coefficients of $X_{10}$ for
$100000 \le n \le 100012$}\label{tab:alphan}
\end{table}
In the notation of \eqref{eqn:gamma}, $\gamma(n) = \alpha(2n) n!/ (2n)!$,
so Table~\ref{tab:alphan} is sufficient to approximate
$J^{d,50000}_{X_{10}, ev}$ for $d \le 6$.  

Let
\begin{align}
A = \mathstrut & 6.288077476637003007403984783011648580 \times 10^{516790} \cr
B = \mathstrut & 1.600352019320098623551973272940701704 \times 10^{21} \cr
C = \mathstrut & 7.155862655552087639602840363255312494 \times 10^{8} .
\end{align}
Then we find that
\begin{equation}
A J^{6,50000}_{X_{10}, ev}(C x + B) =
   -120 + 5.368 x + 180.045 x^2 - 0.894 x^3 - 30.0058 x^4 + 0 x^5 + x^6 .
\end{equation}
For comparison, $H_6(x/2) = -120+ 180 x^2 - 30 x^4 + x^6$, so we see that
the coefficients are close.  

As another example, consider
$\sinc(x) = \sin(x)/x$.  The Taylor coefficients of that function are
easy to compute analytically, making it possible to explore high derivatives.
For the 10-millionth derivative, with $\mathcal A$,
$\mathcal B$, and $\mathcal C$ chosen appropriately, we find
\begin{equation}
\mathcal A J^{6,5000000}_{\sinc, ev}(\mathcal C x + \mathcal B) =
   -120 - 0.536 x + 180.00045 x^2 + 0.089 x^3 - 30.000058 x^4 + 0 x^5 + x^6 ,
\end{equation}
which is even closer to $H_6(x/2)$, and also suggests that the rate of
convergence is on the scale of~$1/\sqrt{n}$.
These examples help support the suggestion that
the appearance of the Hermite polynomials in the even Jensen polynomials is
a universal phenomenon.

\section{On the random matrix conjectures for (derivatives of) L-functions}
We end by addressing the claims that the appearance of the Hermite
polynomials in the even Jensen polynomials has implications for the random
matrix statistics of L-functions.

The first issue is that the Hermite polynomials appear in the
even Jensen polynomials for $\xi^{(n)}(\frac12 + z)$, in other words, in the
classical Jensen polynomials for~$\xi^{(n)}(\frac12 + \sqrt{\mathstrut x})$.
Since the zeros of $\xi^{(n)}(\frac12 + \sqrt{\mathstrut x})$ lie in the
left half-plane, close to the negative real axis, as $n$ increases those
zeros move onto the negative real axis and shift to the left.
In the limit, the zeros fall off the negative real
edge of the complex plane, and suitably rescaled (but not shifted), the limiting function
is~$e^z$.  In the scaled but unshifted classical Jensen polynomials,
all the zeros accumulate at $z=-1$.
Since it is $J^{d,n}_{\xi(\frac12 + \sqrt{\cdot}), cl}(z/d)$
which converges to $\xi^{(n)}(\frac12 + \sqrt{z})$ as $d\to\infty$,
those zeros do not reveal anything about the limiting function
at $z$ in a compact subset.

The second issue is the density of zeros.  The main claim for a connection to
random matrix statistics was that both the zeros of Hermite polynomials and
the eigenvalues of matrices in the Gaussian Unitary Ensemble~(GUE) 
have a density given by the semicircular law.  That is true, but when
using the GUE to model zeros of L-functions, the semicircle density is
a defect, not a feature.  One must artificially rescale the eigenvalues
of GUE matrices to achieve the flat density of zeros of L-functions.
It is the statistics of the spacings, not the density of zeros, which are modeled by
random matrices.  For $x$ in a bounded interval, the zeros of the
Hermite polynomial $H_d(x)$ approach equal spacing as $d\to\infty$.
If those were modeling zeros of derivatives, it would merely be
a reflection of the limiting cosine function, where all information
about the original distribution of zeros has been lost. 


\end{document}